\documentclass[12pt,reqno]{amsart}

\usepackage{amscd,amssymb}

\newcommand{\Z}{{\mathbb Z}}

\newcommand{\C}{{\mathbb C}}
\newcommand{\R}{{\mathbb R}}

\newcommand{\CC}{{\mathcal C}}

\newcommand{\OO}{{\mathcal O}}
\newcommand{\PP}{{\mathcal P}}

\newcommand{\TT}{{\mathcal T}}

\newcommand{\www}{\widetilde}

\newcommand{\paa}{\partial}

\newcommand{\nnn}{\nabla}

\DeclareMathOperator{\rank}{rank}
\DeclareMathOperator{\rk}{rk}

\DeclareMathOperator{\supp}{supp}

\theoremstyle{plain}
\newtheorem{lemma}{Lemma}[section]
\newtheorem{theorem}[lemma]{Theorem}

\newtheorem{corollary}[lemma]{Corollary}

\theoremstyle{definition}
\newtheorem{definition}[lemma]{Definition}

\newtheorem{remark}[lemma]{Remark}
\newtheorem{remarks}[lemma]{Remarks}

\newtheorem{notations}[lemma]{Notations}

\begin{document}

\title[Potentials of a Frobenius like structure]
{Potentials of a Frobenius like structure} 

\author[C. Hertling and A. Varchenko]{Claus Hertling and Alexander Varchenko}

\address{Lehrstuhl f\"ur Mathematik VI,
Universit\"at Mannheim, A5,6, 68131 Mannheim, Germany}
\email{hertling\char64 math.uni-mannheim.de}

\address{Department of Mathematics, University of North Carolina 
at Chapel Hill, Chapel Hill, NC 27599-3250, USA}

\email{anv@email.unc.edu}




\thanks{This work was supported by the DFG grant He2287/4-1
(SISYPH), the second author was supported in part by NSF grant DMS-1362924 and Simons Foundation grant \#336826}

\maketitle


\begin{abstract}
This paper proves the existence of potentials of the first
and second kind of a Frobenius like structure in a frame
which encompasses families of arrangements.

The frame uses the notion of matroids.
For the proof of the existence of the potentials, 
a power series ansatz is made.
The proof that it works requires that certain
decompositions of tuples of coordinate vector fields
are related by certain elementary transformations.
This is shown with a nontrivial result on matroid partition.
\end{abstract}




\section{Introduction and main results}\label{s1}
\setcounter{equation}{0}

\noindent
A Frobenius manifold comes equipped locally with a potential.
If one gives a definition which does not mention this potential
explicitly, one nevertheless obtains it immediately by the following 
elementary fact: Let $z_i$ be the coordinates on $\C^n$
and $\paa_i=\frac{\paa}{\paa z_i}$ be the coordinate vector fields.
Let $M$ be a convex open subset of $\C^n$
and $\TT_M$ be the holomorphic tangent bundle of $M$.
Let $A:\TT_M^3\to\OO_M$ be a symmetric map such that also
$\paa_i A(\paa_j,\paa_k,\paa_l)$ is symmetric in $i,j,k,l$.
Then a potential $F\in\OO_M$ with 
$\paa_i\paa_j\paa_k F=A(\paa_i,\paa_j,\paa_k)$ exists. On Frobenius manifolds see \cite{D, M}.

This paper is devoted to a nontrivial generalization of this fact.
The generalization turns up in the theory of families of arrangements as in \cite[ch. 3]{V2}. 
The geometry there looks at first view similar to the geometry
of Frobenius manifolds, but at second view, it is quite
different.

At first view, one finds in both cases data
$(M,K,\nnn^K,C,S,\zeta)$ with the following properties.
$M$ is an open subset of $\C^n$ (with coordinates $z_i$
and coordinate vector fields $\paa_i=\frac{\paa}{\paa z_i}$).
$K\to M$ is a holomorphic vector bundle with a flat holomorphic
connection $\nnn^K$. $C$ is a Higgs field, i.e. an
$\OO_M$-linear map
\begin{eqnarray}\label{1.1}
C:\OO(K)\to \Omega_M^1\otimes \OO(K)
\end{eqnarray}
such that all the endomorphisms $C_X:K\to K,$ $X\in\TT_M$,
commute: $C_XC_Y=C_YC_X$. And $C$ and $\nnn^K$ satisfy
the integrability condition 
\begin{eqnarray}\label{1.2}
\nnn^K_{\paa_i}C_{\paa_j} = \nnn^K_{\paa_j}C_{\paa_j}
\qquad \textup{for all }i,j\in\{1,...,n\}
\end{eqnarray}
(which is equivalent to $\nnn^K(C)=0$, see remark \ref{t4.1}).
$S$ is a $\nnn^K$-flat symmetric nondegenerate and
Higgs field invariant pairing. $\zeta$ is a global
nowhere vanishing section of $K$.

At second view, one sees the differences.
In the case of a Frobenius manifold, $M$ is the Frobenius
manifold, $\rk K=n$, and (much stronger)
$C_\bullet\zeta:\TT_M\to \OO(K)$ is an isomorphism
and all the sections $C_{\paa_i}\zeta$ are $\nnn^K$-flat.
One obtains an identification of $TM$ with $K$ and of the
coordinate vector fields $\paa_i$ with the flat sections
$C_{\paa_i}\zeta$.

In the case of a family of arrangements, $\rk K\geq n$,
and the $\nnn^K$-flat sections in $K$ have the following
much more surprising form.
Define $J:=\{1,...,n\}$. 
A family of arrangements in $\C^k$ with $k<n$ as in \cite[ch. 3]{V2}
comes equipped with vectors $(v_i)_{i\in J}$ in
$M(1\times k,\C)=\{\textup{row vectors of length }k\textup{ with values in }\C\}$ 
such that $\langle v_1,..,v_n\rangle
=M(1\times k,\C)$.  A subset $\{i_1,...,i_k\}\subset J$
is called {\it maximal independent} if $v_{i_1},...,v_{i_k}$
is a basis of $M(1\times k,\C)$. The sections
$C_{\paa_{i_1}}...C_{\paa_{i_k}}\zeta$ in $K$ 
for such subsets $\{i_1,...,i_k\}$ are $\nnn^K$-flat.

The purpose of this paper is to show that also in this situation
a potential exists which resembles the potential of a
Frobenius manifold.
This is nontrivial. The proof combines the integrability 
condition \eqref{1.2} with intricate combinatorial considerations
which are due to the complicated form of the $\nnn^K$-flat
sections.

Theorem \ref{t1.2} is the main result.
Definition \ref{t1.1} gives the frame and the used notions.
The frame is in two mild aspects more general than the 
data above in the case of arrangements.
First, $S$ is more general, and second, the maximal
independent subsets $\{i_1,...,i_k\}\subset J$ are 
maximal independent with respect to an arbitrary matroid
$(J,F)$ of rank $k$. See definition \ref{t2.1}
for the notion of a matroid.

\begin{definition}\label{t1.1}
(a) A {\it Frobenius like structure of order} 
$(n,k,m)\in\Z_{>0}^3$ 
with $n\geq k$ is a tuple $(M,K,\nnn^K,C,S,\zeta,(J,F))$ 
with the following properties.
$M,K,\nnn^K,C,\zeta$ and $J$ are as above.
$S$ is a $\nnn^K$-flat $m$-linear form $S:\OO(K)^m\to\OO_M$,
which is Higgs field invariant, i.e. 
\begin{eqnarray}\label{1.3}
S(C_Xs_1,s_2,...,s_m)=S(s_1,C_Xs_2,...,s_m)=...=
S(s_1,s_2,...,C_Xs_m)
\end{eqnarray}
for $s_1,s_2,...,s_m\in\OO(K)$ and $X\in \TT_M$.
$(J,F)$ is a matroid with rank $r(J)=k$.
For any maximal independent subset $\{i_1,...,i_k\}\subset J$ 
the section $C_{\paa_{i_1}}...C_{\paa_{i_k}}\zeta$
is $\nnn^K$-flat.

\medskip
(b) Some notations: For any subset $I=\{i_1,...,i_k\}\subset J$, 
the differential operator $\paa_I:=\paa_{i_1}...\paa_{i_k}$ 
and the endomorphism 
$C_I:=C_{\paa_{i_1}}...C_{\paa_{i_k}}:\OO(K)\to\OO(K)$ 
are well defined (they do not depend on the chosen order 
of the elements $i_1,...,i_k$).

\medskip
(c) In the situation of (a), a {\it potential of the first kind}
is a function $Q\in \OO_M$ with
\begin{eqnarray}\label{1.4}
\paa_{I_1}...\paa_{I_m}Q=S(C_{I_1}\zeta,...,C_{I_m}\zeta)
\end{eqnarray}
for any $m$ maximal independent subsets $I_1,...,I_m\subset J$.
A {\it potential of the second kind} is a function 
$L\in\OO_M$ with 
\begin{eqnarray}\label{1.5}
\paa_i\paa_{I_1}...\paa_{I_m}L=S(C_{\paa_i}C_{I_1}\zeta,...,C_{I_m}\zeta)
\end{eqnarray}
for any $m$ maximal independent subsets $I_1,...,I_m\subset J$ 
and any $i\in J$.
\end{definition}

\begin{theorem}\label{t1.2}
Let $(M,K,\nnn^K,C,S,\zeta,(J,F))$ be a Frobenius like 
structure of some order $(n,k,m)\in\Z_{>0}^3$. 
Then locally (i.e. near any $z\in M\subset \C^n$)
potentials of the first and second kind exist.
\end{theorem}

Notice that by formulas (\ref{1.4}) and (\ref{1.5}) the potential of the first kind determines the matrix elements
of the $m$-linear form $S$ on the flat sections $C_{\paa_{i_1}}...C_{\paa_{i_k}}\zeta$
 and the potential of the second kind determines the matrix elements
of the Higgs operators $C_{\paa_{i}}$ acting on the flat sections 
$C_{\paa_{i_1}}...C_{\paa_{i_k}}\zeta$. Thus all information on the $m$-linear form and
the Higgs operators is packed into the two potential functions.

At the end of the paper, several remarks discuss the 
case of arrangements and the relation to Frobenius manifolds. 
In the case of arrangements, one has an
$(n,k,2)$-Frobenius type structure, but also other ingredients,
which lead to richer geometry. 
In the case of a Frobenius manifold, one has an 
$(n,1,2)$-Frobenius type structure. The potential
$L$ above generalizes the potential of a Frobenius manifold.
For generic arrangements,
a global explicit construction of the potentials $Q$ and $L$ 
had been given in \cite{V3}. Recently this was generalized
in \cite{PV} to all families of arrangements as in 
\cite[ch. 3]{V2}.

Section \ref{s2} cites a nontrivial result of J. Edmonds
\cite[4. Theorem]{E} on matroid partition
and adds some considerations. Section \ref{s3}
applies an implication of it to a combinatorial
situation which in turn is needed in the proof of the
main theorem \ref{t1.2} in section \ref{s4}.
Section \ref{s4} concludes with some remarks.

We thank a referee of an earlier version \cite{HV} 
of this paper for pointing us to the result on matroid partition.
This led to the present version of the paper which uses matroids.
The second author thanks  MPI in Bonn for hospitality during his visit in 2015-2016.

\section{Matroid partition}\label{s2}
\setcounter{equation}{0}

\begin{definition}\label{t2.1} (E.g. \cite{E}) 
A {\it matroid} $(E,F)$ is a finite set $E$ together with
a nonempty family $F\subset \PP(E)$ of subsets of $E$,
called {\it independent sets}, such that the following holds.
\begin{list}{}{}
\item[(i)]
Every subset of an independent set is independent.
\item[(ii)]
For every subset $A\subset E$, all maximal independent
subsets of $A$ have the same cardinality, called the
rank $r(A)$ of $A$.
\end{list}
\end{definition}

For example, if $V$ is a vector space and $(v_e)_{e\in E}$
is a tuple of elements which generates $V$,
one obtains a matroid where a subset $B\subset E$ is 
independent if and only if the tuple $(v_b)_{b\in B}$
is a linearly independent tuple of vectors. 
In the case of a family of arrangements, such a matroid will
be used.

The following result on matroid partition 
was proved by J. Edmonds \cite{E}.

\begin{theorem}\label{t2.2} \cite[4. Theorem]{E}.
Let $(E,F_i)$, $i=1,...,m$, be matroids which are defined
on the same set $E$. Let $r_i(A)$ be the rank of $A\subset E$
relative to $(E,F_i)$. The following two conditions
are equivalent.
\begin{list}{}{}
\item[($\alpha$)]
The set $E$ can be partitioned into a family 
$\{I_i\}_{i=1,...,m}$ of sets $I_i\in F_i$.
\item[($\beta$)]
Any set $A\subset E$ satisfies
\begin{eqnarray}\label{2.1}
|A|\leq \sum_{i=1}^m r_i(A).
\end{eqnarray}
\end{list}
\end{theorem}

The implication 
$(\alpha)\Rightarrow(\beta)$ is immediate: Suppose that
$\{I_i\}_{i=1,...,m}$ is a partition of $E$ with
$I_i\in F_i$. Then for any $A\subset E$
\begin{eqnarray*}
A=\dot\bigcup_{i=1}^m A\cap I_i,\quad 
|A|=\sum_{i=1}^m |A\cap I_i|\leq \sum_{i=1}^m r_i(A).
\end{eqnarray*}
But the implication $(\beta)\Rightarrow(\alpha)$ is nontrivial.
The proof in \cite{E} is an involved inductive algorithm.

We are interested in the more special situation in 
theorem \ref{t2.6}.
Before, two lemmata are needed.

\begin{definition}\label{t2.3} \cite{E}
(a) A minimal dependent set of elements of a matroid is called
a {\it circuit}.

(b) For any number $l\in \Z_{\geq 0}$ and any finite set $E$ 
with $|E|\geq l$, the set 
$F^{(l,E)}:=\{I\subset E\, |\, |I|\leq l\}$
defines obviously a matroid $(E,F^{(l,E)})$, the 
{\it uniform matroid} of rank $l$.
\end{definition}

\begin{lemma}\label{t2.4} \cite[Lemma 2]{E}
The union of any independent set $I$ and any element $e$
of a matroid contains at most one circuit of the matroid.
\end{lemma}

\begin{lemma}\label{t2.5}
Let $(E,F)$ be a matroid. 
Let $A_1,A_2\subset E$ be subsets. 
For $i=1,2$, let $I_i\subset A_i$ be a maximal
independent subset of $A_i$. Suppose that $I_1\cup I_2$
is an independent set. Then $I_1\cup I_2$ is a maximal
independent subset of $A_1\cup A_2$, and $I_1\cap I_2$
is a maximal independent subset of $A_1\cap A_2$.
\end{lemma}

{\bf Proof:}
Suppose that for some element $b\in (A_1\cup A_2)-(I_1\cup I_2)$
the union $I_1\cup I_2\cup \{b\}$ is independent.
Then for some $i\in\{1,2\}$, $b\in A_i$.
But $I_i\cup\{b\}$ is a larger independent subset of 
$A_i$ than $I_i$, a contradiction. This proves that
$I_1\cup I_2$ is a maximal independent subset of $A_1\cup A_2$.

Suppose that for some element $b\in (A_1\cap A_2)-(I_1\cap I_2)$
the union $(I_1\cap I_2)\cup\{b\}$ is independent.
If $b\in I_i$ then $b\notin I_j$ where $\{i,j\}=\{1,2\}$.
Then $I_j\cup\{b\}$ is an independent subset of $A_j$,
a contradiction to the maximality of $I_j$. Therefore $b\notin I_1\cup I_2$.
Thus for $i=1,2$, the set $I_i\cup\{b\}\subset A_i$ is dependent
as it is larger than $I_i$. Therefore it contains a circuit
$C_i\subset I_i\cup\{b\}$. Obviously 
$C_i\cap (I_i-I_j)\neq\emptyset$ where $\{i,j\}=\{1,2\}$.
Thus $C_1\neq C_2$. Both are circuits in $(I_1\cup I_2)\cup\{b\}$,
a contradiction to lemma \ref{t2.4}.
This proves that $I_1\cap I_2$ is a maximal independent 
subset of $A_1\cap A_2$.
\hfill$\Box$

\begin{theorem}\label{t2.6}
Let $(E,F_i)$, $i=1,...,m$, be matroids which are defined
on the same set $E$ and which satisfy together 
$(\alpha)$ and $(\beta)$ in theorem \ref{t2.2}. 
Suppose that $F_m=F^{(l,E)}$ for 
some $l\in \Z_{\geq 0}$ with $l\leq |E|$.
Suppose that the set
\begin{eqnarray}\label{2.2}
G:=\{A\subset E\, |\ |A|=l+\sum_{i=1}^{m-1} r_i(A)\}
\end{eqnarray}
contains the set $E$. 

(a) Then this set $G$ is closed under the 
operations union and intersection of sets.
Especially, it contains a set called $A_{min}\subset E$
which is the unique minimal element of $G$ with respect
to the partial order given by inclusion.
Of course $A_{min}\neq\emptyset$ if and only if $l\geq 1$.

(b) Now suppose $l\geq 1$. 
Then $A_{min}=A_{par}$ where $A_{par}$ is the set 
\begin{eqnarray}\label{2.3}
A_{par}&:=&\{b\in E\, |\, \exists\  \textup{a partition }
\{I_i\}_{i=1,...,m}\textup{ of }E\\
&& \hspace*{2cm} 
\textup{such that }I_i\in F_i\textup{ and }b\in I_m\}.
\nonumber
\end{eqnarray}
\end{theorem}

{\bf Proof:} (a) 
Choose a partition $\{I_i\}_{i=1,...,m}$ of $E$ with 
$I_i\in F_i$. 
For any subset $A\subset E$, it induces a partition
$A=\dot\bigcup_{i=1}^m A\cap I_i$ of $A$ into
subsets $(A\cap I_i)\in F_i$. 
If $A\in G$, then by \eqref{2.2} each set $A\cap I_i$
is a maximal independent subset of $A$ with respect
to the matroid $(E,F_i)$. 
As $|A|\geq l$, especially $|A\cap I_m|=l$. 
As $E$ itself is in $G$, $|I_m|=l$, and thus
$A\cap I_m =I_m$ for any set $A\in G$.

Let $A_1,A_2\in G$. 
For any $i=1,...,m$, lemma \ref{t2.5}
applies to the maximal independent sets $A_1\cap I_i$
and $A_2\cap I_i$ of $A_1$ respectively $A_2$
relative to the matroid $(E,F_i)$,
because also $(A_1\cup A_2)\cap I_i\in F_i$. 
Therefore $(A_1\cup A_2)\cap I_i$ is a maximal independent
subset of $A_1\cup A_2$ relative to $(E,F_i)$,
and $(A_1\cap A_2)\cap I_i$ is a maximal independent
subset of $A_1\cap A_2$ relative to $(E,F_i)$. 
Also, $I_m=A_1\cap I_m=A_2\cap I_m$ shows
\begin{eqnarray*}
I_m=(A_1\cup A_2)\cap I_m =(A_1\cap A_2)\cap I_m.
\end{eqnarray*}
Now $A_1\cup A_2\in G$ and $A_1\cap A_2\in G$ are obvious.
Therefore $G$ is closed under the operations union
and intersection of sets. 

(b) $A_{par}\subset A_{min}$: 
Fix an arbitrary element $b\in A_{par}$.
Choose a partition $\{I_i\}_{i=1,...,m}$ of $E$ with 
$I_i\in F_i$ and $b\in I_m$. 
Recall $A_{min}\cap I_m=I_m$. Thus $b\in A_{min}$.

$A_{min}\subset A_{par}$: 
Fix an arbitrary element $b\in A_{min}$. 
Define $\www E:=E-\{b\}$.
Any set $A\subset\www E$ does not
contain $A_{min}$, because $b\in A_{min}$. 
Therefore any set $A\subset\www E$ satisfies $A\notin G$ and 
\begin{eqnarray}\label{2.4}
|A|\leq -1+l+\sum_{i=1}^{m-1}r_i(A).
\end{eqnarray}
Consider the matroids $(\www E,\www F_i)$,
where $\www F_i:=\{I\in F_i\, |\, b\notin I\}$
for $i\in\{1,...,m-1\}$ and $\www F_m:=F^{(l-1,\www E)}$. 
For $i\in\{1,...,m-1\}$ the rank of $A\subset\www E$ 
relative to $(\www E,\www F_i)$ is equal to the rank $r_i(A)$
of $A$ relative to $(E,F_i)$.

By \eqref{2.4} and theorem \ref{t2.2}, 
a partition $\{\www I_i\}_{i=1,...,m}$
of $\www E$ with $\www I_i\in \www F_i$ exists.
Now the sets $I_i:=\www I_i$ for $i=1,...,m-1$,
and $I_m:=\www I_m\cup\{b\}$ form a partition of $E$
with $I_i\in F_i$. This shows $b\in A_{par}$.
\hfill$\Box$

\section{An equivalence between index systems}\label{s3}
\setcounter{equation}{0}

\noindent 
In this section we fix three positive integers 
$n,k,m\in\Z_{>0}$ with $n\geq k$ and a matroid $(J,F)$ 
with underlying set $J=\{1,...,n\}$, 
rank function $r:\PP(J)\to\Z_{\geq 0}$ and rank $r(J)=k$.

\begin{notations}\label{t3.1}
As usual $\Z^J:=\{\textup{maps}:J\to\Z\}$
and $\Z_{\geq 0}^J:=\{\textup{maps}:J\to\Z_{\geq 0}\}$.
The set $\Z^J$ is an additive group,
the set $\Z_{\geq 0}^J$ is an additive monoid.

For $j\in J$ denote by $[j]\in\Z_{\geq 0}^J$
the map with $[j](j)=1$ and $[j](i)=0$ for any $i\neq j$.
Then any map $T\in \Z^J$ can be written as
$T=\sum_{j=1}^nT(j)\cdot [j]$. 
For $T\in\Z^J$ denote $|T|:=\sum_{j=1}^nT(j)\in\Z$. 
The support of $T\in \Z^J$ is 
$\supp T:=\{j\in J\, |\, T(j)\neq 0\}$. 
The map
\begin{eqnarray}\label{3.1}
d_H:\Z^J\times\Z^J\to\Z_{\geq 0},\quad (T_1,T_2)\mapsto
\sum_{j\in J}|T_1(j)-T_2(j)|
\end{eqnarray}
is a metric on $\Z^J$. 
On $\Z^J$ one has the partial ordering $\leq$ with
\begin{eqnarray}\label{3.2}
S\leq T\iff S(j)\leq T(j)\quad\forall\ j\in J.
\end{eqnarray}
Any map $T\in \Z_{\geq 0}^J$ with $|T|=t\in\Z_{\geq 0}$
is called a {\it system of elements of $J$} 
or simply a {\it system} or a {\it $t$-system}.
If $S$ and $T$ are systems with $S\leq T$, then
$S$ is a {\it subsystem} of $T$.
\end{notations}

\begin{definition}\label{t3.2} 
Here $l\in\Z_{\geq 0}$. Here all systems are systems
of elements of $J$.
\begin{list}{}{}
\item[(a)] 
A system $T\in\Z_{\geq 0}^J$ is a {\it base} if $\supp T\in F$ 
and  $|T|=k$ (so the support $\supp T$ is a maximal independent
subset of $J$ and all $T(a)\in\{0;1\}$).

\item[(b)]
A {\it strong decomposition} of an $(mk+l)$-system $T$ is a 
decomposition $T=T^{(1)}+...+T^{(m+1)}$ into $m$ $k$-systems 
$T^{(1)},...,T^{(m)}$ and one $l$-system $T^{(m+1)}$ 
such that $T^{(1)},...,T^{(m)}$ are
bases (and $T^{(m+1)}$ is an arbitrary $l$-system; 
e.g. if $l=0$ then $T^{(m+1)}=0$ automatically).

\item[(c)] 
An $(mk+l)$-system is {\it strong} if it admits a strong decomposition.

\item[(d)]
A {\it good decomposition} of an $N$-system $T$ 
with $N\geq mk+1$ is a decomposition $T=T_1+T_2$ 
into two systems such that $T_2$ is a strong $(mk+1)$-system 
of elements of $J$.

\item[(e)]
Two good decompositions $T_1+T_2=T$ and $S_1+S_2=T$ 
of an $N$-system $T$ with $N\geq mk+1$
are {\it locally related}, 
notation: $(S_1,S_2)\sim_{loc} (T_1,T_2)$, 
if there are strong decompositions 
$S^{(1)}_2+...+S^{(m+1)}_2=S_2$ of $S_2$ and 
$T^{(1)}_2+...+T^{(m+1)}_2=T_2$ of $T_2$ with 
$S^{(j)}_2=T^{(j)}_2$ for $1\leq j\leq m$. 
Of course, $\sim_{loc}$ is a reflexive and symmetric relation.

\item[(f)]
Two good decompositions $T_1+T_2=T$ and $S_1+S_2=T$ 
of an $N$-system $T$ with $N\geq mk+1$
are {\it equivalent}, 
notation: $(S_1,S_2)\sim (T_1,T_2)$, 
if there is a sequence $\sigma_1,\sigma_2,...,\sigma_r$
for some $r\in\Z_{\geq 1}$ of good decompositions of $T$
such that $\sigma_1=(S_1,S_2)$, $\sigma_r=(T_1,T_2)$ and
$\sigma_j\sim_{loc}\sigma_{j+1}$ for $j=1,...,r-1$.
Of course, $\sim$ is an equivalence relation.

\end{list}
\end{definition}

The main result of this section is the following theorem \ref{t3.3}.

\begin{theorem}\label{t3.3}
Let $T\in\Z_{\geq 0}^J$ be an $N$-system for some $N\geq mk+1$
which has good decompositions. Then all its good decompositions
are equivalent.
\end{theorem}

The theorem will be proved after the proofs
of corollary \ref{t3.4} and lemma \ref{t3.5}.
Corollary \ref{t3.4} is a corollary of theorem \ref{t2.6}.

\begin{corollary}\label{t3.4}
Fix a strong $(mk+l)$-system $T\in\Z_{\geq 0}^J$ 
with $l\in\Z_{\geq 0}$. 
Then for any $B\subset J$
\begin{eqnarray}\label{3.3}
\sum_{j\in B}T(j)\leq l+m\cdot r(B).
\end{eqnarray}
The set
\begin{eqnarray}\label{3.4}
G(T)&:=& \{B\subset \supp T\, |\, \sum_{j\in B}T(j)
= l+m \cdot r(B)\}
\end{eqnarray}
contains $\supp T$ and is closed under the operations union and 
intersection of sets. Especially, it contains a set 
called $A_{min}(T)\subset \supp T$ which is the unique
minimal element with respect to inclusion.
In the case $l\geq 1$, define the set
\begin{eqnarray}\label{3.5}
A_{dec}(T)&:=&\{b\in J\, |\, \exists\ 
\textup{a strong decomposition }\\
&&T=T^{(1)}+...+T^{(m+1)}\textup{ with }b\in \supp T^{(m+1)}\} .\nonumber
\end{eqnarray}
Then $A_{min}(T)=A_{dec}(T)$.
\end{corollary}

{\bf Proof:}  
We will construct from $T$ certain lifts of the
matroids $(J,F)$ and $(J,F^{(l,J)})$ to matroids on
the set $E:=\{1,2,...,mk+l\}$ and go with them into
theorem \ref{t2.6}. 
Choose a map $f:E\to J$ with $|f^{-1}(j)|=T(j)$.
Define the sets 
\begin{eqnarray*}
F_1=...=F_m&:=& \{A\subset E\, |\, f|_A:A\to J
\textup{ injective,}\ f(A)\in F\}\subset\PP(E),  \\
F_{m+1}&:=& F^{(l,E)}\subset\PP(E).
\end{eqnarray*}
Then $(E,F_i)$ for $i\in\{1,...,m+1\}$ is a matroid.
Together they satisfy $(\alpha)$ in theorem \ref{t2.2}
(with $m+1$ instead of $m$)
because $T$ is a strong $(mk+l)$-system.
We go into theorem \ref{t2.6} with $m+1$ instead of $m$.

That $T$ is a strong $(mk+l)$-system, gives also
$E\in G$ and \eqref{3.3}.

Therefore the set $A_{min}$ in theorem \ref{t2.6} is 
well defined. The set $A_{par}$ is well defined,
anyway. One sees easily 
\begin{eqnarray*}
r_1(A)=...=r_m(A)&=&r(f(A))\quad\textup{for }A\subset E,\\
G&=&\{f^{-1}(B)\, |\, B\in G(T)\}.
\end{eqnarray*}
Therefore $G(T)$ contains $\supp T$ and is closed under
the operations union and intersection of sets.
Now one sees also easily
\begin{eqnarray*}
A_{min}&=&f^{-1}(A_{min}(T)),\quad 
A_{par}=f^{-1}(A_{dec}(T)),
\end{eqnarray*}
and thus $A_{min}(T)=A_{dec}(T)$.
\hfill$\Box$

\begin{lemma}\label{t3.5}
Let $S$ and $T\in \Z_{\geq 0}^J$ be two strong $(mk+1)$-systems.
At least one of the following two alternatives holds.

\begin{list}{}{}
\item[$(\alpha)$] 
$T$ has a strong decomposition $T=T^{(1)}+...+T^{(m+1)}$
with $T^{(m+1)}=[i]$ for some $i\in \supp T$ with $T(i)>S(i)$.
\item[$(\beta)$] 
For any strong decomposition $S=S^{(1)}+...+S^{(m+1)}$
a strong decomposition $T=T^{(1)}+...+T^{(m+1)}$ 
with $T^{(m+1)}=S^{(m+1)}$ exists.
\end{list}
\end{lemma}

{\bf Proof:}
Suppose that $(\alpha)$ does not hold.
Then for any $i\in A_{dec}(T)$ $S(i)\geq T(i)$. Especially
\begin{eqnarray*}
\sum_{i\in A_{dec}(T)}S(i)&\geq&  \sum_{i\in A_{dec}(T)}T(i)
= 1+m\cdot r(A_{dec}(T)).
\end{eqnarray*}
The equality uses $A_{dec}(T)=A_{min}(T)\in G(T)$.
Now \eqref{3.3} for $S$ instead of $T$ shows that $\geq$ can be replaced by $=$.
Therefore $A_{dec}(T)\in G(S)$. 
Any element of $G(S)$ contains $A_{min}(S)$. This and
the equality $A_{dec}(S)=A_{min}(S)$ give
$$A_{dec}(S)=A_{min}(S)\subset A_{dec}(T).$$
Thus $(\beta)$ holds.
\hfill$\Box$

\bigskip
{\bf Proof of theorem \ref{t3.3}:}
Let $(S_1,S_2)$ and $(T_1,T_2)$ be two different good decompositions
of an $N$-system $T$ of elements of $J$ (with $N\geq mk+1$).
Then $S_2$ and $T_2$ are strong $(mk+1)$-systems of elements 
of $J$.
At least one of the two alternatives $(\alpha)$ and $(\beta)$ in
lemma \ref{t3.5} holds for $S_2$ and $T_2$.

\medskip
{\bf First case, $(\alpha)$ holds:} 
Let $T_2=T_2^{(1)}+...+T_2^{(m+1)}$ be a strong decomposition
with $T_2^{(m+1)}=[i]$ for some $i\in\supp T_2$ with 
$T_2(i)>S_2(i)$.
Then a $j\in \supp T$ with $T_1(j)>S_1(j)$ and $T_2(j)<S_2(j)$ exists. 
The decomposition 
\begin{eqnarray}\label{3.6}
T=R_1+R_2\quad\textup{with }R_1=T_1-[j]+[i],\quad R_2=T_2+[j]-[i]
\end{eqnarray}
is a good decomposition of $T$ because 
$T_2^{(1)}+...+T_2^{(m)}+[j]$ is a strong decomposition of $R_2$. 
The good decompositions $(R_1,R_2)$ and $(T_1,T_2)$ are locally
related, $(R_1,R_2)\sim_{loc}(T_1,T_2)$, and thus equivalent,
\begin{eqnarray}\label{3.7}
(R_1,R_2)\sim(T_1,T_2).
\end{eqnarray}
Furthermore, 
\begin{eqnarray}\label{3.8}
d_H(R_2,S_2)=d_H(T_2,S_2)-2.
\end{eqnarray}

\medskip
{\bf Second case, $(\beta)$ holds:}
Let $T_2=T_2^{(1)}+...+T_2^{(m+1)}$ and 
$S_2=S_2^{(1)}+...+S_2^{(m+1)}$
be strong decompositions of $T_2$ and $S_2$ with 
$T_2^{(m+1)}=S_2^{(m+1)}=[a]$ for some $a\in \supp T$.
Two elements $b,c\in\supp T$ with 
$T_1(b)>S_1(b)$, $T_2(b)<S_2(b)$, and $T_1(c)<S_1(c)$, $T_2(c)>S_2(c)$
exist. 
Consider the decompositions of $T$ and $S$,
\begin{eqnarray}\label{3.9}
T&=&R_1+R_2\quad\textup{with }R_1=T_1-[b]+[a],R_2=T_2+[b]-[a],\\
S&=&Q_1+Q_2\quad\textup{with }Q_1=S_1-[c]+[a],Q_2=S_2+[c]-[a].\label{3.10}
\end{eqnarray}
They are good decompositions because $R_2$ has the strong
decomposition $R_2=T^{(1)}+...+T^{(m)}+[b]$ 
and $Q_2$ has the strong decomposition
$Q_2=S^{(1)}+...+S^{(m)}+[c]$. 
The local relations
\begin{eqnarray*}
(R_1,R_2)\sim_{loc} (T_1,T_2)\quad\textup{and}\quad 
(Q_1,Q_2)\sim_{loc}(S_1,S_2)
\end{eqnarray*}
and the equivalences
\begin{eqnarray}\label{3.11}
(R_1,R_2)\sim (T_1,T_2)\quad\textup{and}\quad 
(Q_1,Q_2)\sim (S_1,S_2)
\end{eqnarray}
hold. Furthermore
\begin{eqnarray}\label{3.12}
d_H(R_2,Q_2)=d_H(T_2,S_2)-2.
\end{eqnarray}

\medskip
The properties \eqref{3.7}, \eqref{3.8}, \eqref{3.11} and \eqref{3.12}
show that in both cases the equivalence classes of $(S_1,S_2)$ and
$(T_1,T_2)$ contain good decompositions whose second members
are closer to one another with respect to the metric $d_H$
than $T_2$ and $S_2$. This shows that $(S_1,S_2)$ and $(T_1,T_2)$
are in one equivalence class.  \hfill$\Box$

\section{Potentials of the first and second kind}\label{s4}
\setcounter{equation}{0}

The main part of this section is devoted to the proof of
theorem \ref{t1.2}. At the end some remarks on the
relation to families of arrangements and 
Frobenius manifolds are made.

\begin{remark}\label{t4.1}
Here a coordinate free formulation of the integrability
condition \eqref{1.2} will be given.
For $M,\nnn^K$ and $C$ as in the introduction,
$\nnn^K(C)\in \Omega^2_M\otimes\OO(\textup{End}(K))$ 
is the 2-form on $M$
with values in $\textup{End}(K)$ such that for $X,Y\in\TT_M$
\begin{eqnarray}\label{4.1}
\nnn^K(C)(X,Y)&=& \nnn^K_X(C_Y)-\nnn^K_Y(C_X)-C_{[X,Y]}.
\end{eqnarray}
Now \eqref{1.2} is equivalent to $\nnn^K(C)=0$
\end{remark}

{\bf Proof of theorem \ref{t1.2}:}
Let $(M,K,\nnn^K,C,S,\zeta,(J,F))$ be a Frobenius 
like structure of some order $(n,k,m)\in\Z_{>0}^3$.

We need some notations. If $T\in\Z_{\geq 0}^J$ is a system
of elements of $J$, then 
\begin{eqnarray*}
(z-x)^T&:=&\prod_{i\in J}(z_i-x_i)^{T(i)}\quad\textup{for any }x\in\C^n,\\
T!:=\prod_{i\in J}T(i)!,\quad
\paa_T&:=&\prod_{i\in J}\paa_{z_i}^{T(i)},
\quad C_T:=\prod_{i\in J}C_{\paa_{z_i}}^{T(i)}.
\end{eqnarray*}
Thus, if $S$ and $T$ are systems of elements of $J$, then
\begin{eqnarray}
\paa_T(z-x)^S=\left\{\begin{array}{ll}
0&\textup{ if }T\not\leq S,\\
\frac{S!}{(S-T)!}\cdot (z-x)^{S-T}& \textup{ if }T\leq S,
\end{array}\right. \label{4.2}
\end{eqnarray}
for any $x\in\C^n$.

\medskip
The existence of a (not just local, but even global) 
potential $Q$ of the first kind is trivial.
The function
\begin{eqnarray}\label{4.3}
Q&:=& \sum_{T\textup{ with }(*)}\frac{1}{T!}\cdot 
S(C_T \zeta,\zeta,...,\zeta)\cdot z^T \quad(m\textup{ times }\zeta),
\hspace*{1cm}\\
(*)&:& T\in\Z_{\geq 0}^J
\textup{ is a strong }mk\textup{-system (definition \ref{t3.1}(c))}.\nonumber
\end{eqnarray}
works. It is a homogeneous polynomial of degree $mk$
and contains only monomials which are relevant for \eqref{1.2}.
In fact, one can add to this $Q$ an arbitrary linear combination
of the monomials $z^T$ for the $mk$-systems $T$ which are not strong,
so which are not relevant for \eqref{1.2}.

\medskip
The existence of a potential $L$ of the second kind is not trivial.
Let some $x\in M$ be given. We make the power series ansatz
\begin{eqnarray}\label{4.4}
L&:=& \sum_{T\in\Z_{\geq 0}^J} a_T\cdot (z-x)^T,
\end{eqnarray}
where the coefficients $a_T$ have to be determined.
If $T$ satisfies $|T|\leq mk$ or if it satisfies
$|T|\geq mk+1$, but does not admit a good decomposition
(definition \ref{t3.1} (d)), then the conditions \eqref{1.3}
are empty for $a_T(z-x)^T$ because of \eqref{4.2},
so then $a_T$ can be chosen arbitrarily, e.g. $a_T:=0$ works.

Now consider $T$ with $|T|\geq mk+1$ which admits good decompositions.
Then each good decomposition $T=T_1+T_2$ gives 
via \eqref{1.3} a candidate
\begin{eqnarray}\label{4.5}
a_T(T_1,T_2)&:=& \frac{1}{T!}\cdot \left(\paa_{T_1}
S(C_{T_2}\zeta,\zeta,...,\zeta)\right)(x),
\end{eqnarray}
for the coefficient $a_T$ of $(z-x)^T$ in $L$.
We have to show that the candidates $a_T(T_1,T_2)$ for all
good decompositions $(T_1,T_2)$ of $T$ coincide.

Suppose that two good decompositions $(T_1,T_2)$ and $(S_1,S_2)$
are locally related, $(T_1,T_2)\sim_{loc}(S_1,S_2)$
(definition \ref{t3.1} (e)), but not equal.
Then there are strong decompositions $T_2=T_2^{(1)}+...+T_2^{(m)}+[a]$
and $S_2=T_2^{(1)}+...+T_2^{(m)}+[b]$ with $a\neq b$,
and thus also $T_1-[b]=S_1-[a]\in\Z_{\geq 0}^J$ holds.
Because any $T_2^{(j)}$, $j\in\{1,...,m\}$, is independent,
$C_{T_2^{(j)}}\zeta$ is $\nnn^K$-flat.
This and \eqref{4.3} give
\begin{eqnarray}
&&\paa_{z_b}S(C_{T_2}\zeta,\zeta,...,\zeta)\nonumber\\
&=& \paa_{z_b}S(C_{\paa_{z_a}}C_{T_2^{(1)}}\zeta, C_{T_2^{(2)}}\zeta,...,
C_{T_2^{(m)}}\zeta)\nonumber\\
&=& S(\nnn^K_{\paa_{z_b}}(C_{\paa_{z_a}})C_{T_2^{(1)}}\zeta, C_{T_2^{(2)}}\zeta,...,
C_{T_2^{(m)}}\zeta)\nonumber\\
&=& S(\nnn^K_{\paa_{z_a}}(C_{\paa_{z_b}})C_{T_2^{(1)}}\zeta, C_{T_2^{(2)}}\zeta,...,
C_{T_2^{(m)}}\zeta)\nonumber\\
&=& \paa_{z_a}S(C_{\paa_{z_b}}C_{T_2^{(1)}}\zeta, C_{T_2^{(2)}}\zeta,...,
C_{T_2^{(m)}}\zeta)\nonumber\\
&=&\paa_{z_a}S(C_{S_2}\zeta,\zeta,...,\zeta).
\label{4.6}
\end{eqnarray}
This implies 
\begin{eqnarray}\label{4.7}
a_T(T_1,T_2) = a_T(S_1,S_2),
\end{eqnarray}
so the locally related good decompositions $(T_1,T_2)$ and $(S_1,S_2)$
give the same candidate for $a_T$. 
Thus all equivalent (definition \ref{t3.1} (f)) good decompositions
give the same candidate for $a_T$.
By theorem \ref{t3.3}, all good decompositions of $T$ are equivalent.
Therefore they all give the same candidate for $a_T$.
Thus a potential $L$ of the second kind exists as a formal power series
as in \eqref{4.4}.

It is in fact a convergent power series because of the following.
There are finitely many strong $(mk+1)$-systems $T_2$. 
Each determines the coefficients $a_T$ for all $T\geq T_2$. 
We put $a_T:=0$ for $T$ which do not admit
good decompositions. The part of $L$ in \eqref{4.4} which is determined by some strong $(mk+1)$-system $T_2$ 
is a convergent power series.
Thus $L$ is the {\it union} of finitely many overlapping
convergent power series. It is easy to see that it is itself convergent. 
This finishes the proof of theorem \ref{t1.2}.
\hfill$\Box$

\begin{remark}\label{t4.2}
In \cite[ch. 3]{V2} families of arrangements are considered which
give rise to Frobenius like structures 
$(M,K,\nnn^K,C,S,\zeta,(J,F))$ of order $(n,k,2)$,
see the special case of generic arrangements in \cite{V1,V3}.

Start with two positive integers $k$ and $n$ with $k<n$
and with a matrix $B:=(b_i^j)_{i=1,..,n;j=1,..,k}\in M(n\times k,\C)$
with $\rank B=k$. Define $J:=\{1,...,n\}$. 
Here the matroid $(J,F)$ is the {\it vector matroid}
(also called {\it linear matroid}) of the tuple
$(v_i)_{i\in J}$ of row vectors 
$v_i:= (b^j_i)_{j=1,...k}$ of the matrix $B$.
More precisely, a subset $A\subset J$ is independent,
if the tuple $(v_i)_{i\in A}$ is a 
linearly independent system of vectors.

Consider $\C^n\times \C^k$ with the
coordinates $(z,t)=(z_1,...,z_n,t_1,...,t_k)$ and with the projection
$\pi:\C^n\times \C^k\to\C^n$. Define the functions
\begin{eqnarray}\label{4.8}
g_i:=\sum_{j=1}^kb_i^j\cdot t_j,\quad f_i:=g_i+z_i
\quad\textup{for }i\in J
\end{eqnarray}
on $\C^n\times \C^k$. 

We obtain on $\C^n\times \C^k$ the arrangement 
$\CC=\{H_i\}_{i\in J}$, where $H_i$ is the zero set of $f_i$.
Let $U(\CC):=\C^n\times \C^k-\bigcup_{i\in J}H_i$ 
be the complement.
For every $x\in\C^n$, the arrangement $\CC$ restricts to an 
arrangement $\CC(x)$ on $\pi^{-1}(x)\cong\C^k$. 
For almost all $x\in\C^n$ the arrangement $\CC(x)$ 
is {\it essential} (definition in \cite{V2})
with normal crossings. The subset $\Delta\subset\C^n$ where
this does not hold, is a hypersurface and is called 
the {\it discriminant},
see \cite[3.2]{V2}. Define $M:=\C^n-\Delta$.

A set $I=\{i_1,...,i_k\}\subset J$ is maximal independent, i.e.
$(v_{i_1},...,v_{i_k})$ is a basis of $M(1\times k,\C)$, 
if and only if for some (or equivalently for any) $x\in \C^n$ 
the hyperplanes $H_{i_1}(x),...,H_{i_k}(x)$ are transversal.

Let $a=(a_1,...,a_n)\in(\C^*)^n$ be a system of {\it weights}
such that for any $x\in M$ the weighted arrangement
$(\CC(x),a)$ is {\it unbalanced}: See \cite{V2} for the definition
of {\it unbalanced}, e.g. $a\in\R_{>0}^n$ is unbalanced,
also a generic system of weights is     unbalanced.
The {\it master function} of the weighted arrangement $(\CC,a)$
is
\begin{eqnarray}\label{4.9}
\Phi_a(z,t):=\sum_{i\in J}a_i\log f_i.
\end{eqnarray}
Several deep facts are related to this master function.
We use some of them in the following. See \cite{V2} for references.

For $z\in M$ all critical points of $\Phi_{a}$ are isolated,
and the sum $\mu$ of their Milnor numbers is independent of 
the unbalanced weight $a$ and the parameter $z\in M$. 
The bundle 
\begin{eqnarray}\label{4.10}
K:=\bigcup_{z\in M}K_z\quad\textup{with }
K_z:=\OO(U(\CC)\cap\pi^{-1}(z))/\left(\frac{\paa\Phi_a}{\paa t_j}\, |\,
j=1,...,k\right)
\end{eqnarray}
over $M$ is a vector bundle of $\mu$-dimensional algebras.

It comes equipped with the section $\zeta$ of unit elements 
$\zeta(z)\in K_z$, a Higgs field $C$, 
a {\it combinatorial connection} $\nnn^K$ and a pairing $S$.
The Higgs field $C:\OO(K)\to \Omega^1_M\otimes \OO(K)$
is defined with the help of the period map
\begin{eqnarray}\label{4.11}
\Psi:TM\to K,\quad \paa_{z_i}\mapsto \left[\frac{\paa\Phi_a}{\paa z_i}\right]
=\left[\frac{a_i}{f_i}\right]=:p_i
\end{eqnarray}
by
\begin{eqnarray}\label{4.12}
C_{\paa_{z_i}}(h):=p_i\cdot h\qquad\textup{ for }h\in K_z.
\end{eqnarray}
Because of 
\begin{eqnarray}\label{4.13}
0=\left[\frac{\paa\Phi_a}{\paa t_j}\right]
=\sum_{i=1}^n b^j_i p_i,
\end{eqnarray}
the Higgs field vanishes on the vector fields
$X_j:=\sum_{i=1}^n b^j_i\paa_i$, $j\in\{1,...,k\}$, 
\begin{eqnarray}\label{4.14}
C_{X_j}=0\qquad\textup{for }j\in\{1,...,k\}.
\end{eqnarray}
In fact the whole geometry of the family of arrangements is 
invariant with respect to the flows of these vector fields.

The sections $\det(b_i^j)_{i\in I,j=1,...,k}\cdot C_I\zeta$ 
for all maximal independent sets $I=\{i_1,...,i_k\}\subset J$
generate the bundle $K$, and they satisfy only relations with 
constant coefficients in $\Z$. 
The combinatorial connection $\nnn^K$ is the unique
flat connection such that the sections $C_I\zeta$ for $I\subset J$
 maximal independent are $\nnn^K$-flat.
The sections $\det(b_i^j)_{i\in I,j=1,...,k}\cdot C_I\zeta$ 
for $I\subset J$ maximal independent generate a 
$\nnn^K$-flat $\Z$-lattice structure on $K$.

The pairing $S$ comes from the Grothendieck residue with respect
to the volume form 
\begin{eqnarray}\label{4.15}
\frac{dt_1\land...\land dt_k}{\prod_{j=1}^k 
\frac{\paa\Phi_a}{\paa t_j}}.
\end{eqnarray}
It is symmetric, nondegenerate, $\nnn^K$-flat, 
multiplication invariant and Higgs field invariant.

\smallskip
The existence of potentials of the first and second kind for families of arrangements was conjectured in \cite{V1}. If all the $k\times k$ minors of the matrix $B=(b_i^j)$ are nonzero, the potentials were constructed in  \cite{V1}, cf. \cite{V3}. 
In \cite{PV} this was generalized to all cases in this remark
\ref{t4.2}. 
The potentials are given by explicit formulas in terms of the linear functions defining the hyperplanes in $\C^n$ composing the discriminant.

\end{remark}

\begin{remarks}\label{t4.3}
(i) The situation in remark \ref{t4.2} is in several aspects
richer than a Frobenius like structure of type $(n,k,m)$.
The bundle $K$ is a bundle of algebras.
The sections $C_I\zeta$ for  maximal independent sets 
$I\subset J$ generate the bundle. The sections  
$\det(b_i^j)_{i\in I,j=1,...,k}\cdot C_I\zeta$
generate a flat $\Z$-lattice structure in $K$.
The Higgs field vanishes on the vector fields $X_1,...,X_k$.
The $m$-linear form $S$ is a pairing ($m=2$) and is nondegenerate.
We will not discuss the $\Z$-lattice structure, but we will discuss
some logical relations between the other enrichments and some 
implications of them.

\medskip
(ii) Let $(M,K,\nnn^K,C,S,\zeta,V,(v_1,...,v_n))$ be a Frobenius like
structure of order $(n,k,m)$. Suppose that it satisfies the 
{\it generation condition}
\begin{eqnarray}\label{4.16}
\text{(GC)}&& \textup{The sections }C_I\zeta
\textup{ for maximal independent sets } I\subset J\\
&&\textup{generate the bundle }K. \nonumber
\end{eqnarray}
Let $\mu$ be the rank of $K$. Then for any $x\in M$, the endomorphisms
$C_X,X\in T_xM$, generate a $\mu$-dimensional commutative subalgebra
$A_z\subset\textup{End}(K_x)$. And any endomorphism which commutes
with them is contained in this subalgebra. This gives a rank $\mu$
bundle $A$ of commutative algebras. And the map
\begin{eqnarray}\label{4.17}
A\to K,\quad B\mapsto B\zeta,
\end{eqnarray}
is an isomorphism of vector bundles and induces a commutative and
associative multiplication on $K_x$ for any $x\in M$, with unit field
$\zeta(x)$. Therefore the special section $\zeta$  and 
the generation condition (GC), which exist and hold in remark 
\ref{t4.2}, give the multiplication on the bundle $K$ there.

\medskip
(iii) In the situation in (ii) with the condition (GC),
the $m$-linear form is multiplication invariant because it is 
Higgs field invariant. The condition (GC) implies also that it is symmetric:
\begin{eqnarray*}
S(C_{I_1}\zeta,C_{I_2}\zeta,...,C_{I_m}\zeta)
=S(C_{I_{\sigma(1)}}\zeta,C_{I_{\sigma(2)}}\zeta,...,C_{I_{\sigma(m)}}\zeta)
\end{eqnarray*}
for any  maximal independent sets $I_1,...,I_m$ 
and any permutation $\sigma\in S_m$.

\medskip
(iv) The following special case gives rise to
Frobenius manifolds without Euler fields.
Consider a Frobenius like structure
$(M,K,\nnn^K,C,S,\zeta,(J,F))$ of order
$(n,1,2)$ with nondegenerate pairing $S$, 
$\nnn^K$-flat section $\zeta$, 
the uniform matroid $(J,F)=(J,F^{(1,J)})$ and the condition
that the map $C_\bullet \zeta:TM\to K$ is an isomorphism.
Then the sections $C_{\paa_i}\zeta$ generate
the bundle $K$ and are $\nnn^K$-flat.
Here $M$ becomes a Frobenius manifold (without Euler field)
whose flat structure is the naive flat structure
of $\C^n\supset M$.
The potential $L$ is the potential of the Frobenius manifold.
\end{remarks}

\end{document}